\documentclass[12pt,reqno]{amsart}

\usepackage[dvips]{graphicx} 

\usepackage{amssymb, latexsym, amsmath, amscd, array,
}

\usepackage{fancyhdr}

\newtheorem{theorem}{Theorem}[section]

\theoremstyle{definition}

\newtheorem{question}[theorem]{Question}

\numberwithin{equation}{section}
\newcommand\N {{\mathbb N}} 

\newcommand\R {{\mathbb R}}
\newcommand\Z {{\mathbb Z}} 
\newcommand\st{{\rm st}} 

\global\advance\count0 by 2

\begin{document}

\thispagestyle{empty}
\rightline{\em TMME, vol.~7, no.~1, p.~\thepage} %first page

\vskip1in

{\centerline{\Huge\bf When is .999... less than 1?}}

\bigskip\bigskip

{\centerline{\large\bf Karin Usadi Katz and Mikhail G. Katz$^{0}$}}

\address{Department of Mathematics, Bar Ilan University, Ramat Gan
52900 Israel} \email{katzmik@macs.biu.ac.il}

\footnotetext{%
\medskip\noindent 
Supported by the Israel Science Foundation (grants no.~84/03 and
1294/06).\hfil\break
{\bf\em The Montana Mathematics Enthusiast, ISSN 1551-3440, Vol. 7,
No. 1, pp.~3--30} \hfil\break
2010 \copyright The Montana Council of Teachers of Mathematics \&
Information Age Publishing
}



\bigskip\bigskip\bigskip\noindent
We examine alternative interpretations of the symbol described as {\em
nought, point, nine recurring\/}.  Is ``an infinite number of 9s''
merely a figure of speech?  How are such alternative interpretations
related to infinite cardinalities?  How are they expressed in
Lightstone's ``semicolon'' notation?  Is it possible to choose a
canonical alternative interpretation?  Should unital evaluation of the
symbol~$.999\ldots$ be inculcated in a pre-limit teaching environment?
The problem of the unital evaluation is hereby examined from the
pre-$\R$, pre-$\lim$ viewpoint of the {\em student\/}.

%\maketitle{}

%\tableofcontents

\section{Introduction}%1

Leading education researcher and mathematician D. Tall \cite{Ta10}
comments that a mathematician ``may think of the physical line as an
{\em approximation\/} to the infinity of numbers, so that the line is
a practical representation of numbers[, and] of the number line as a
visual representation of a precise numerical system of decimals.''
Tall concludes that ``this still does not alter the fact that there
are connections in the minds of students based on experiences with the
number line that differ from the formal theory of real numbers and
cause them to feel confused.''

One specific experience has proved particularly confusing to the
students, namely their encounter with the evaluation of the
symbol~$.999\ldots$ to the standard real value~$1$.  Such an
evaluation will be henceforth referred to as the {\em unital
evaluation\/}.

We have argued \cite{KK} that the students are being needlessly
confused by a premature emphasis on the unital evaluation, and that
their persistent intuition that~$.999\ldots$ can fall short of~$1$,
can be rigorously justified.  Other interpretations (than the unital
evaluation) of the symbol~$.999\ldots$ are possible, that are more in
line with the students' naive initial intuition, persistently reported
by teachers.  From this viewpoint, attempts to inculcate the
equality~$.999\ldots=1$ in a teaching environment prior to the
introduction of limits (as proposed in~\cite{WAD}), appear to be
premature.

To be sure, certain student intuitions are clearly dysfunctional, such
as a perception that $\forall\epsilon\exists\delta$ and
$\exists\delta\forall\epsilon$ are basically ``the same thing''.  Such
intuitions need to be uprooted.  However, a student who intuits
$.999\ldots$ as a dynamic process (see R.~Ely \cite{El, El10}) that
never quite reaches its final address, is grappling with a fruitful
cognitive issue at the level of the world of proceptual symbolism, see
Tall \cite{Ta10}.  The student's functional intuition can be channeled
toward mastering a higher level of abstraction at a later, limits/$\R$
stage.

In '72, A. Harold Lightstone published a text entitled {\em
Infinitesimals\/} in the {\em American Mathematical Monthly}
\cite{Li}.  If~$\epsilon>0$ is infinitesimal (see Appendix~\ref{ten}),
then~$1-\epsilon$ is less than~$1$, and Lightstone's extended decimal
expansion of~$1- \epsilon$ starts with more than any finite number of
repeated~$9$s.%
\footnote{\label{avigad}
By the work of J.~Avigad \cite{Av01}, the phenomenon can already be
expressed in primitive recursive arithmetic, in the context of
Skolem's non-standard models of arithmetic, see answer to
Questions~\ref{q61} and \ref{skolem} below.}
Such a phenomenon was briefly mentioned by Sad, Teixeira, and Baldino
\cite[p.~286]{STB}.

The symbol~$.999\ldots$ is often said to possess more than any finite
number of~$9$s.  However, describing the real decimal~$.999\ldots$ as
possessing an {\em infinite number of~$9$s\/} is only a figure of
speech, as {\em infinity\/} is not a number in standard analysis.

The above comments have provoked a series of thoughtful questions from
colleagues, as illustrated below in a question and answer format
perfected by Imre Lakatos in \cite{La2}.

\section{Frequently 
asked questions: when is~$.999\ldots$ less than~$1$?}%2

\begin{question}%1.1
Aren't there many standard proofs that $0.999\ldots =1$?  Since we
can't have that and also~$0.999\ldots \not=1$ at the same time, if
mathematics is consistent, then isn't there necessarily a flaw in your
proof?
\end{question}

Answer. The standard proofs%
\footnote{The proof exploiting the long division of $1$ by $3$, is
dealt with in the answer to Question~\ref{q35}.}
are of course correct, in the context of the standard real numbers.
However, the annals of the teaching of decimal notation are full of
evidence of student frustration with the unital evaluation
of~$.999\ldots$ This does not mean that we should tell the students an
untruth.  What this does mean is that it may be instructive to examine
why exactly the students are frustrated.  

\begin{question}%1.1
Why are the students frustrated?
\end{question}

Answer.  The important observation here is that the students are not
told about either of the following two items:
\begin{enumerate}
\item
the real number system; 
\item
limits,
\end{enumerate}
before they are exposed to the topic of decimal notation, as well as
the problem of unital evaluation.  What we point out is that so long
as the number system has not been specified explicitly, the students'
hunch that~$.999\ldots$ falls infinitesimally short of 1 can be
justified in a rigorous fashion, in the framework of Abraham
Robinson's \cite{Ro66, Ro} non-standard analysis.$^{\ref{avigad}}$

%\begin{figure}%[ht]
%\includegraphics[height=3in]{robinson.ps}
%\caption{Robinson on non-standard integers in~$\N^* \setminus \N$}
%\label{robinsonfig}
%\end{figure}

\begin{question}%1.2
Isn't it a problem with the proof that the
definitions aren't precise?  You say that~$0.999\ldots$ has an
``unbounded number of repeated digits~$9$".  That is not a meaningful
mathematical statement; there is no such number as ``unbounded".  If
it is to be precise, then you need to provide a formal definition of
``unbounded", which you have not done.
\end{question}

Answer.  The comment was not meant to be a precise definition.  The
precise definition of~$.999\ldots$ as a real number is well known.  The
thrust of the argument is that {\em before} the number system has been
explicitly specified, one can reasonably consider that the ellipsis
``..." in the symbol~$.999\ldots$ is in fact ambiguous.  From this point
of view, the notation~$.999\ldots$ stands, not for a single number, but
for a class of numbers,%
\footnote{For a more specific choice of such a number, see the answer
to Question~\ref{q36}.}
all but one of which are less than~$1$.

Note that F. Richman \cite{Ri} creates a natural semiring (in the
context of decimal expansions), motivated by constructivist
considerations, where certain cancellations are disallowed (as they
involve infinite ``carry-over'').  The absence of certain
cancellations (i.e. subtractions) leads to a system where a strict
inequality~$.999\ldots < 1$ is satisfied.  The advantage of the
hyperreal approach is that the number system remains a field, together
with the extension principle and the transfer principle (see
Appendix~\ref{ten}).

\begin{question}%1.3
Doesn't decimal representation have the same meaning in standard
analysis as non-standard analysis?
\end{question}

Answer.  Yes and no.  Lightstone \cite{Li} has developed an extended
decimal notation%
\footnote{For details see Appendix~\ref{ten}, item~\ref{1011} below.}
that gives more precise information about the hyperreal.  In his
notation, the standard real $.999\ldots$ would appear as
\[
.999\ldots;\ldots 999\ldots
\]

\begin{question}%1.4
Since non-standard analysis is a conservative extension of the
standard reals, shouldn't all existing properties of the standard
reals continue to hold?
\end{question}

Answer.  Certainly, $.999\ldots;\ldots 999 \ldots$ equals $1$, on the
nose, in the hyperreal number system, as well.  An accessible account
of the hyperreals can be found in chapter~6: {\em Ghosts of departed
quantities\/} of Ian Stewart's popular book {\em From here to
infinity\/} \cite{St}.  In his unique way, Stewart has captured the
essense of the issue as follows in \cite[p.~176]{St09}:
\begin{quote}
The standard analysis answer is to take `$\ldots$' as indicating
passage to a limit.  But in non-standard analysis there are many
different interpretations.
\end{quote}
In particular, a {\em terminating\/} infinite
decimal~$.999\ldots;\ldots 999$ is less than~$1$.

\begin{question}%1.5
Your expression ``terminating infinite decimals'' sounds like
gibberish.  How many decimal places do they have exactly?  How can
infinity terminate?
\end{question}

Answer.  If you are troubled by this, you are in good company.
A remarkable passage by Leibniz is a testimony to the enduring appeal
of the metaphor of infinity, even in its, paradoxically, {\em
terminated\/} form.  In a letter to Johann Bernoulli dating from june
1698 (as quoted in Jesseph \cite[Section~5]{Je}), Leibniz speculated
concerning
\begin{quote}
lines [...] which are terminated at either end, but which nevertheless
are to our ordinary lines, as an infinite to a finite.  
\end{quote}
He further speculates as to the possibility of
\begin{quote}
a point in space which can not be reached in an assignable time by
uniform motion.  And it will similarly be required to conceive a time
terminated on both sides, which nevertheless is infinite, and even
that there can be given a certain kind of eternity [...] which is
terminated.
\end{quote}
Ultimately, Leibniz rejected any metaphysical reality of such
quantities, and conceived of both infinitesimals and infinite
quantities as ideal numbers, falling short of the reality of the
familiar appreciable quantities.%
\footnote{Ely \cite{El10} presents a case study of a student who
naturally developed an intuitive system of infinitesimals and
infinitely large quantities, bearing a striking resemblance to
Leibniz's system.  Ely concludes: ``By recognizing that some student
conceptions that appear to be misconceptions are in fact nonstandard
conceptions, we can see meaningful connections between cognitive
structures and mathematical structures of the present and past that
otherwise would have been overlooked.''}

\begin{question}%1.5
You say that in Lightstone's notation, the nonstandard number
represented by $.999\ldots;\ldots 99\hat{9}$ is less~$1$.  Wouldn't he
consider this as something different from~$.999\ldots$, since he uses
a different notation, and that he would say
\[
.999\ldots;\ldots 99\hat{9} < .999\ldots  = 1?
\]
\end{question}

Answer. Certainly.

\begin{question}%1.6
Aren't you arbitrarily redefining~$.999\ldots$ as equal to the
non-standard number~$.999\ldots;\ldots 99\hat{9}$, which would
contradict the standard definition?
\end{question}

Answer.  No, the contention is that the ellipsis notation is
ambiguous, particularly as perceived by pre-$\lim$, pre-$\R$ students.
The notation could reasonably be applied to a class of numbers%
\footnote{For a more specific choice of such a number, see the answer
to Question~\ref{q36}.}
infinitely close to~$1$.

\begin{question}%1.7
You claim that ``there is an unbounded number of 9s in~$.999\ldots$,
but saying that it has infinitely many $9$s is only a figure of
speech".  Now there are several problems with such a claim.  First,
there is no such object as an ``unbounded number".  Second,
``infinitely many 9s" not a figure of speech, but rather quite
precise.  Doesn't ``infinite'' in this context mean the countable
cardinal number,~$\aleph_0$ in Cantor's notation?
\end{question}

Answer.  One can certainly choose to call the output of a series
whatever one wishes.  The terminology ``infinite sum" is a useful and
intuitive term, when it comes to understanding standard calculus.  In
other ways, it can be misleading.  Thus, the term contains no hint of
the fact that such an ``$\aleph_0$-fold sum" is only a partial
operation, unlike the inductively defined~$n$-fold sums.  Namely, a
series can diverge, in which case the infinite sum is undefined (to be
sure, this does not happen for decimal series representing real
numbers).  

Furthermore, the ``$\aleph_0$-fold sum" intuition creates an
impediment to understanding Lightstone extended decimals
\[
 .a_1 a_2 a_3 \ldots ; \ldots a_H \ldots
\]
If one thinks of the standard real as an~$\aleph_0$-fold sum of the
countably many terms such as~$a_1/10, a_2/100, a_3/1000$, etc., then
it may appear as though Lightstone's extended decimals {\em add\/}
additional positive (infinitesimal) terms to the real value one
started with (which seems to be already ``present" to the left of the
semicolon).  It then becomes difficult to understand how such an
extended decimal can represent a number {\em less\/} than~$1$.

For this reason, it becomes necessary to analyze the {\em infinite
sum\/} figure of speech, with an emphasis on the built-in {\em limit}.

\begin{question}%1.8
Are you trying to convince me that the expression {\em infinite sum},
routinely used in Calculus, is only a figure of speech?
\end{question}

Answer.  The debate over whether or not an {\em infinite sum\/} is a
figure of speech, is in a way a re-enactment of the foundational
debates at the end of the 17th and the beginning of the 18th century,
generally thought of as a Newton-Berkeley debate.%
\footnote{
See footnote~\ref{bernoulli} for a historical clarification.}
The founders of the calculus thought of
\begin{enumerate}
\item
the derivative as a ratio of a pair of infinitesimals, and of 
\item
the integral as an infinite sum of terms~$f(x)dx$.
\end{enumerate}
Bishop Berkeley \cite{Be} most famously%
\footnote{
\label{rolle}Similar criticisms were expressed by Rolle,
thirty years earlier; see Schubring~\cite{Sch}.}
criticized the former in terms of the familiar {\em ghosts of departed
quantities\/} (see \cite[Chapter~6]{St}) as follows.  The
infinitesimal
\[
dx
\]
appearing in the denominator is expected, at the beginning of the
calculation, to be nonzero (the {\em ghosts\/}), yet at the end of the
calculation it is neglected as if it were zero (hence, {\em
departed\/} quantities).  The implied stripping away of an
infinitesimal at the end of the calculation occurs in evaluating an
integral, as well.

To summarize, the integral is not an infinite Riemann sum, but rather
the standard part of the latter (see Section~\ref{ten},
item~\ref{1012}).  From this viewpoint, calling it an infinite sum is
merely a figure of speech, as the crucial, final step is left out.

A.~Robinson solved the 300-year-old logical inconsistency of the
infinitesimal definition of the integral, in terms of the standard
part function.%
\footnote{
See Appendix~\ref{ten}, item~\ref{103}.}

\begin{question}%
Hasn't historian Bos criticized Robinson for being excessive in
enlisting Leibniz for his cause?
\end{question}

Answer.  In his essay on Leibniz, H. Bos \cite[p.~13]{Bos}
acknowledged that Robinson's hyperreals provide
\begin{quote}
[a] preliminary explanation of why the calculus could develop on the
insecure foundation of the acceptance of infinitely small and
infinitely large quantities.
\end{quote}
F.~Medvedev \cite{Me87, Me98} further points out that
\begin{quote}
[n]onstandard analysis makes it possible to answer a delicate question
bound up with earlier approaches to the history of classical analysis.
If infinitely small and infinitely large magnitudes are regarded as
inconsistent notions, how could they [have] serve[d] as a basis for
the construction of so [magnificent] an edifice of one of the most
important mathematical disciplines?
\end{quote}

\section{Is {\em infinite sum\/} a figure of speech?}%3

\begin{question}%
Perhaps the historical definition of an integral, as an infinite sum
of infinitesimals, had been a figure of speech.  But why is an
infinite sum of a sequence of real numbers more of a figure of speech
than a sum of two real numbers?
\end{question}

Answer.  Foundationally speaking, the two issues (integral and
infinite sum) are closely related.  The series can be described
cognitively as a proceptual encapsulation of a dynamic process
suggested by a sequence of finite sums, see D.~Tall \cite{Ta10}.
Mathematically speaking, the convergence of a series relies on the
completeness of the reals, a result whose difficulty is of an entirely
different order compared to what is typically offered as arguments in
favor of unital evaluation.

The rigorous justification of the notion of an integral is identical
to the rigorous justification of the notion of a series.  One can
accomplish it finitistically with epsilontics, or one can accomplish
it infinitesimally with standard part.  In either case, one is dealing
with an issue of an entirely different nature, as compared to
finite~$n$-fold sums.

\begin{question}%2.2
You have claimed that ``saying that it has an infinite number of 9s is
only a figure of speech".  Of course ``infinity" is not a number in
standard analysis: this word refers to a number in the cardinal number
system, i.e. the cardinality%
\footnote{
See more on cardinals in answer to Question~\ref{q61}.}
of the number of digits; it does not refer to a number in the real
number system.
\end{question}

Answer.  One can certainly consider an infinite string of 9s labeled
by the standard natural numbers.  However, when challenged to write
down a precise definition of $.999\ldots$, one invariably falls back
upon the {\em limit\/} concept (and presumably the respectable
epsilon, delta definition thereof).  Thus, it turns out that
$.999\ldots$ is really the limit of the sequence~$.9$,~$.99$,~$.999$,
etc.%
\footnote{
The related sequence $.3$,~$.33$,~$.333$, etc.~is discussed
in the answer to Question~\ref{q35}.}
Note that such a definition never uses an infinite string of $9$s
labeled by the standard natural numbers, but only finite fragments
thereof.  

Informally, when the students are confronted with the problem of the
unital evaluation, they are told that the decimal in question is zero,
point, followed by infinitely many 9s.  Well, taken literally, this
describes the hyperreal number
\[
.999\ldots;\ldots 999000
\]
perfectly well: we have zero, point, followed by~$H$-infinitely many
9s.  Moreover, this statement in a way is truer than the one about the
standard decimal, as explained above (the infinite string is never
used in the actual standard definition).  The hyperreal {\em is\/} an
infinite sum, on the nose.  It is not a limit of finite sums.

\begin{question}%2.3
Do limits have a role in the hyperreal approach?
\end{question}

Answer.  Certainly.  Let~$u_1=.9, u_2=.99, u_3=.999$, etc.  Then the
limit, from the hyperreal viewpoint, is the standard part of~$u_H$ for
any infinite hyperinteger~$H$.  The standard part strips away the
(negative) infinitesimal, resulting in the standard value~$1$, and the
students are right almost everywhere.

\begin{question}%2.4
A mathematical notation is whatever it is defined to be, no more and
no less.  Isn't~$.999\ldots$ defined to be equal to~$1$?
\end{question}

Answer.  As far as teaching is concerned, it is not necessarily up to
research mathematicians to decide what is good notation and what is
not, but rather should be determined by the teaching profession and
its needs, particularly when it comes to students who have not yet
been introduced to~$\R$ and~$\lim$.

\begin{question}%2.5
In its normal context,~$.999\ldots$ is defined unambiguously, shouldn't
it therefore be taught as a single mathematical object?
\end{question}

Answer.  Indeed, in the context of ZFC standard reals and the
appropriate notion of limit, the definition is unambiguous.  The issue
here is elsewhere: what does~$.999\ldots$ {\em look like\/} to
highschoolers when they are exposed to the problem of unital
evaluation, {\em before\/} learning about~$\R$ and~$\lim$.%
\footnote{
A preferred choice of a hyperreal evaluation of the symbol
``$.999\ldots$'' is described in the answers to Questions~\ref{q45},
\ref{q36}, and \ref{q46}.}

\begin{question}%2.6
Don't standard analysis texts provide a unique definition of
$.999\ldots$ that is almost universally accepted, as a certain
infinite sum that (independently) happens to evaluate to 1?
\end{question}

Answer.  More precisely, it is a limit of finite sums, whereas
``infinite sum" is a figurative way of describing the limit.  Note
that the hyperreal sum from 1 to~$H$, where~$H$ is an infinite
hyperinteger, can also be described as an ``infinite sum'', or more
precisely $H$-infinite sum, for a choice of a hypernatural number $H$.

\begin{question}%3.7
There are certain operations that happen to work with ``formal"
manipulation, such as dividing each digit by 3 to result%
\footnote{
The long division of $1$ by $3$ and its implications for unital
evaluation are discussed in detail in the answer to
Question~\ref{q35}.}
in~$0.333\ldots$ But shouldn't such manipulation be taught as
merely a convenient shortcut that happens to work but needs to be
verified independently with a rigorous argument before it is accepted?
\end{question}

Answer.  Correct.  The best rigorous argument, of course, is that the
sequence 
\[
\langle .9, .99, .999, \ldots\rangle
\]
gets closer and closer to $1$ (and therefore~$1$ is the limit by
definition).  The students would most likely find the remark before
the parenthesis, unobjectionable.  Meanwhile, the parenthetical remark
is unintelligible to them, unless they have already taken calculus.

\section{Meanings, standard and non-standard}%4

\begin{question}%4.1
Isn't it very misleading to change the standard meaning
of~$.999\ldots$, even though it may be convenient?  This is in the
context of standard analysis, since non-standard analysis is not
taught very often because it has its own set of issues and
complexities.
\end{question}

Answer.  In the fall of '08, a course in calculus was taught using
H. Jerome Keisler's textbook {\em Elementary Calculus\/} \cite{Ke}.
The course was taught to a group of 25 freshmen.  The TA had to be
trained as well, as the material was new to the TA.  The students have
never been so excited about learning calculus, according to repeated
reports from the TA.  Two of the students happened to be highschool
teachers (they were somewhat exceptional in an otherwise teenage
class).  They said they were so excited about the new approach that
they had already started using infinitesimals in teaching basic
calculus to their 12th graders.  After the class was over, the TA paid
a special visit to the professor's office, so as to place a request
that the following year, the course should be taught using the same
approach.  Furthermore, the TA volunteered to speak to the chairman
personally, so as to impress upon him the advantages of non-standard
calculus.  The~$.999\ldots$ issue was not emphasized in the class.%
\footnote{
Hyperreal pedagogy is analyzed in the answer to Question~\ref{q28}.}

\begin{question}%4.2
Non-standard calculus?  Didn't Errett Bishop explain already that
non-standard calculus constituted a debasement of meaning?
\end{question}

Answer.  Bishop did refer to non-standard calculus as a {\em
debasement of meaning\/} in his {\em Crisis\/} text \cite{Bi75} from
'75.  He clarified what it was exactly he had in mind when he used
this expression, in his {\em Schizophrenia\/} text~\cite{Bi85}.  The
latter text was distributed two years earlier, more precisely in '73,
according to M. Rosenblatt~\cite[p.~ix]{Ros}.  Bishop writes as
follows \cite[p.~1]{Bi85}:
\begin{quote}
Brouwer's criticisms of {\bf classical mathematics} [emphasis
added--MK] were concerned with what I shall refer to as ``the
debasement of meaning''.
\end{quote}
In Bishop's own words, the {\em debasement of meaning\/} expression,
employed in his {\em Crisis\/} text to refer to non-standard calculus,
was initially launched as a criticism of {\em classical mathematics\/}
as a whole.  Thus his criticism of non-standard calculus was
foundationally, not pedagogically, motivated.  

In a way, Bishop is criticizing apples for not being oranges: the
critic (Bishop) and the criticized (Robinson's non-standard analysis)
do not share a common foundational framework.  Bishop's preoccupation
with the extirpation of the law of excluded middle (LEM)%
\footnote{
A defining feature of both intuitionism and Bishop's
constructivism is a rejection of LEM; see footnote~\ref{semimystical}
for a discussion of Bishop's foundational posture within the spectrum
of Intuitionistic sensibilities.}
led him to criticize classical mathematics as a whole in as vitriolic%
\footnote{
\label{vit}Historian of mathematics J.~Dauben noted the vitriolic
nature of Bishop's remarks, see~\cite[p.~139]{Da96}; M.~Artigue
\cite{Ar} described them as {\em virulent\/}; D.~Tall \cite{Ta01}, as
{\em extreme\/}.}
a manner as his criticism of non-standard analysis.

\begin{question}%2.11a
Something here does not add up.  If Bishop was opposed to the rest of
classical mathematics, as well, why did he reserve special vitriol for
his book review of Keisler's textbook on non-standard calculus?
\end{question}

Answer.  Non-standard analysis presents a formidable philosophical
challenge to Bishopian constructivism, which may, in fact, have been
anticipated by Bishop himself in his foundational speculations, as we
explain below.

While Bishop's constructive mathematics (unlike Brouwer's
intuitionism%
\footnote{
Bishop rejected both Kronecker's finitism and Brouwer's
theory of the continuum.}
) is uniquely concerned with finite operations on the
integers, Bishop himself has speculated that ``the primacy of the
integers is not absolute'' \cite[p.~53]{Bi68}:
\begin{quote}
It is an {\bf empirical fact} [emphasis added--MK] that all [finitely
performable abstract calculations] reduce to operations with the
integers.  There is no reason mathematics should not concern itself
with finitely performable abstract operations of other kinds, in the
event that such are ever discovered [...]
\end{quote}
Bishop hereby acknowledges that the {\em primacy of the integers\/} is
merely an {\em empirical fact\/}, i.e.~an empirical observation, with
the implication that the observation could be contradicted by novel
mathematical developments.  Non-standard analysis, and particularly
non-standard calculus, may have been one such development.

\begin{question}%2.10
How is a theory of infinitesimals such a novel development?
\end{question}

Answer.  Perhaps Bishop sensed that a rigorous theory of
infinitesimals is both
\begin{itemize}
\item
not reducible to finite calculations on the integers, and yet
\item
accomodates a finite performance of abstract operations,
\end{itemize}
thereby satisfying his requirements for coherent mathematics.  Having
made a foundational commitment to the {\em primacy of the integers\/}
(a state of mind known as {\em integr-ity\/} in Bishopian
constructivism; see \cite[p.~4]{Bi85}) through his own work and that
of his disciples starting in the late sixties, Bishop may have found
it quite impossible, in the seventies, to acknowledge the existence of
``finitely performable abstract operations of other kinds''.

Birkhoff reports that Bishop's talk at the workshop was not
well-received.%
\footnote{
See Dauben~\cite[p.~133]{Da96} in the name of Birkhoff
\cite[p.~505]{Pr75}.}
The list of people who challenged him (on a number of points) in the
question-and-answer session that followed the talk, looks like the
who-is-who of 20th century mathematics.

\begin{question}%2.10
Why didn't all those luminaries challenge Bishop's {\em debasement\/}
of non-standard analysis?
\end{question}

Answer. The reason is a startling one: there was, in fact, nothing to
challenge him on.  Bishop did not say a word about non-standard
analysis in his oral presentation, according to a workshop
participant~\cite{Man} who attended his talk.%
\footnote{
The participant in question, historian of mathematics K.~Manning, was
expecting just this sort of critical comment about non-standard
analysis from Bishop, but the comment never came.  Manning wrote as
follows on the subject of Bishop's statement on non-standard calculus
published in the written version \cite{Bi75} of his talk: ``I do not
remember that any such statement was made at the workshop and doubt
seriously that it was in fact made.  I would have pursued the issue
vigorously, since I had a particular point of view about the
introduction of non-standard analysis into calculus.  I had been
considering that question somewhat in my attempts to understand
various standards of rigor in mathematics.  The statement would have
fired me up.''}
Bishop appears to have added the {\em debasement\/} comment after the
workshop, at the galley proof stage of publication.  This helps
explain the absence of any critical reaction to such {\em
debasement\/} on the part of the audience in the discussion session,
included at the end of the published version of Bishop's talk.

\begin{question}%4.6
On what grounds did Bishop criticize classical mathematics as
deficient in numerical meaning?
\end{question}

Answer. The quest for greater numerical meaning is a compelling
objective for many mathematicians.  Thus, as an alternative to an
indirect proof (relying on LEM) of the irrationality of $\sqrt{2}$,
one may favor a direct proof of a concrete lower bound, such
as~$\frac{1}{3n^2}$ for the error $|\sqrt{2} - \frac{m}{n}|$ involved.
Bishop discusses this example in \cite[p.~18]{Bi85}.  More generally,
one can develop a methodology that seeks to enhance classical
arguments by eliminating the reliance on LEM, with an attendant
increase in numerical meaning.  Such a methodology can be a useful
{\em companion\/} to classical mathematics.

\begin{question}%4.7
Given such commendable goals, why haven't mainstream mathematicians
adopted Bishop's constructivism?
\end{question}

Answer.  The problem starts when LEM-extirpation is elevated to the
status of the supreme good, regardless of whether it is to the
benefit, or detriment, of numerical meaning.  Such a radical, anti-LEM
species of constructivism tends to be posited, not as a {\em
companion\/}, but as an {\em alternative\/}, to classical mathematics.
Philosopher of mathematics G.~Hellman \cite[p. 222]{Hel93a} notes that
``some of Bishop's remarks (1967) suggest that his position belongs in
[the radical constructivist] category''.  

For instance, Bishop wrote \cite[p.~54]{Bi68} that ``[v]ery possibly
classical mathematics will cease to exist as an independent
discipline.''  He challenged his precursor Brouwer himself, by
describing the latter's theory of the continuum as a ``semimystical
theory'' \cite[p.~10]{Bi67}.  Bishop went as far as evoking the term
``schizophrenia in contemporary mathematics'', see \cite{Bi85}.

\begin{question}
How can the elimination of the law of excluded middle be detrimental
to numerical meaning?
\end{question}

Answer.  In the context of a discussion of the differentiation
procedures in Leibniz's infinitesimal calculus, D.~Jesseph
\cite[Section~1]{Je} points out that
\begin{quote}
[t]he algorithmic character of this procedure is especially important,
for it makes the calculus applicable to a vast array of curves whose
study had previously been undertaken in a piecemeal fashion, without
an underlying unity of approach.
\end{quote}
The algorithmic, computational, numerical meaning of such computations
persists after infinitesimals are made rigorous in Robinson's
approach, relying as it does on classical logic, incorporating the law
of excluded middle (LEM).

To cite an additional example, note that Euclid himself has recently
been found lacking, constructively speaking, by M.~Beeson
\cite{Bee09}.  The latter rewrote as many of Euclid's geometric
constructions as he could while avoiding ``test-for-equality''
constructions (which rely on LEM).  What is the status of those
results of Euclid that resisted Beeson's reconstructivisation?  Are we
prepared to reject Euclid's constructions as lacking in meaning, or
are we, rather, to conclude that their meaning is of a post-LEM kind?

\begin{question}%4.9
Are there examples of post-LEM numerical meaning from contemporary
research?
\end{question}

Answer.  In contemporary proof theory, the technique of proof mining
is due to Kohlenbach, see \cite{KO}.  A logical analysis of classical
proofs (i.e., proofs relying on classical logic) by means of a
proof-theoretic technique known as proof mining, yields explicit
numerical bounds for rates of convergence, see also Avigad \cite{Av}.

\begin{question}%4.10
But hasn't Bishop shown that meaningful mathematics is mathematics
done constructively?
\end{question}

Answer.  If he did, it was by a sleight-of-hand of a successive
reduction of the meaning of ``meaning''.  First, ``meaning'' in a
lofty epistemological sense is reduced to ``numerical meaning''.  Then
``numerical meaning'' is further reduced to the avoidance of LEM.

\begin{question}
Haven't Brouwer and Bishop criticized formalism for stripping
mathematics of any meaning?
\end{question}

Answer. Thinking of formalism in such terms is a common misconception.
The fallacy was carefully analyzed by Avigad and Reck~\cite{AR}.  

From the cognitive point of view, the gist of the matter was
summarized in an accessible fashion by D. Tall
\cite[chapter~12]{Ta10}:
\begin{quote}
The aim of a formal approach is not the stripping away of all human
intuition to give absolute proof, but the careful organisation of
formal techniques to support human creativity and build ever more
powerful systems of mathematical thinking.
\end{quote}
Hilbert sought to provide a finitistic foundation for mathematical
activity, at the meta-mathematical level.  He was prompted to seek
such a foundation as an alternative to set theory, due to the famous
paradoxes of set theory, with ``the ghost of Kronecker''
(see~\cite{AR}) a constant concern.  Hilbert's finitism was, in part,
a way of answering Kronecker's concerns (which, with hindsight, can be
described as intuitionistic/constructive).

Hilbert's program does not entail any denial of meaning at the
mathematical level.  A striking example mentioned by S.~Novikov
\cite{No2} is Hilbert's Lagrangian for general relativity, a deep and
{\em meaningful\/} contribution to both mathematics and physics.
Unfortunately, excessive rhetoric in the heat of debate against
Brouwer had given rise to the famous quotes, which do not truly
represent Hilbert's position, as argued in \cite{AR}.  Hilbert's
Lagrangian may in the end be Hilbert's most potent criticism of
Brouwer, as variational principles in physics as yet have no
intuitionistic framework, see \cite[p.~22]{Bee}.

\begin{question}%4.12
Why would one want to complicate the students' lives by introducing
infinitesimals?  Aren't the real numbers complicated enough?
\end{question}

Answer. The traditional approach to calculus using Weierstrassian
epsilontics (the epsilon-delta approach) is a formidable challenge to
even the gifted students.%
\footnote{
including Paul Halmos; see \cite{Sf}, as well as the answer
to Question~\ref{51} below.}
Infinitesimals provide a means of simplifying the technical aspect of
calculus, so that more time can be devoted to conceptual issues.

\section{Halmos on  infinitesimal subtleties}%5

\begin{question}%5.1
\label{51}
Aren't you exaggerating the difficulty of Weierstrassian epsilontics,
as you call it?  If it is so hard, why hasn't the mathematical
community discovered this until now?
\end{question}

Answer. Your assumption is incorrect.  Some of our best and brightest
have not only acknowledged the difficulty of teaching Weierstrassian
epsilontics, but have gone as far as admitting their own difficulty in
learning it!  For example, Paul Halmos recalls in his autobiography
\cite[p.~47]{Ha}:
\begin{quote}
... I was a student, sometimes pretty good and sometimes less good.
Symbols didn't bother me.  I could juggle them quite well ...[but] I
was stumped by the infinitesimal subtleties of epsilonic analysis.  I
could read analytic proofs, remember them if I made an effort, and
reproduce them, sort of, but I didn't really know what was going on.
\end{quote}
(quoted in A.~Sfard \cite[p.~44]{Sf}).  The eventual resolution of
such pangs in Halmos' case is documented by Albers and Alexanderson
\cite[p.~123]{AA}:
\begin{quote}
... one afternoon something happened ...  suddenly I understood
epsilon.  I understood what limits were ...  All of that stuff that
previously had not made any sense became obvious ...
\end{quote}
Is Halmos' liberating experience shared by a majority of the students
of Weierstrassian epsilontics?

\begin{question}%2.11bb
I don't know, but how can one possibly present a construction of the
hyperreals to the students?
\end{question}

Answer.  You are surely aware of the fact that the construction of the
reals (Cauchy sequences or Dedekind cuts) is not presented in a
typical standard calculus class.%
\footnote{
The issue of constructing number systems is discussed further in the
answers to Questions~\ref{q51} and \ref{q33}.}
Rather, the instructor relies on intuitive descriptions, judging
correctly that there is no reason to get bogged down in
technicalities.  There is no more reason to present a construction of
infinitesimals, either, so long as the students are given clear ideas
as to how to perform arithmetic operations on infinitesimals, finite
numbers, and infinite numbers.  This replaces the rules for
manipulating limits found in the standard approach.%
\footnote{
See the answer to Question~\ref{q33} for more details on the
ultrapower construction.}

\begin{question}%
Non-standard analysis?  Didn't Halmos explain already that it is too
special?
\end{question} 

Answer.  P.~Halmos did describe non-standard analysis as {\em a
special tool, too special\/} \cite[p.~204]{Ha}.  In fact, his
anxiousness to evaluate Robinson's theory may have involved a conflict
of interests.  In the early '60s, Bernstein and Robinson~\cite{BR}
developed a non-standard proof of an important case of the invariant
subspace conjecture of Halmos', and sent him a preprint.  In a race
against time, Halmos produced a standard translation of the
Bernstein-Robinson argument, in time for the translation to appear in
the same issue of {\em Pacific Journal of Mathematics\/}, alongside
the original.  Halmos invested considerable emotional energy (and {\em
sweat\/}, as he memorably puts it in his autobiography%
\footnote{
Halmos wrote \cite[p.~204]{Ha}: ``The Bernstein-Robinson proof [of the
invariant subspace conjecture of Halmos'] uses non-standard models of
higher order predicate languages, and when [Robinson] sent me his
reprint I really had to sweat to pinpoint and translate its
mathematical insight.''}
) into his translation.  Whether or not he was capable of subsequently
maintaining enough of a detached distance in order to formulate an
unbiased evaluation of non-standard analysis, his blunt unflattering
comments appear to retroactively justify his translationist attempt to
deflect the impact of one of the first spectacular applications of
Robinson's theory.

\begin{question}%
How would one express the number~$\pi$ in Lightstone's
``.999\ldots;\ldots 999" notation?
\end{question}

Answer.  Certainly, as follows:
\[
3.141\ldots;\ldots d_{H-1} d_H d_{H+1} \ldots
\]
The digits of a standard real appearing after the semicolon are, to a
considerable extent, determined by the digits before the semicolon.
The following interesting fact might begin to clarify the situation.
Let
\[
d_{\min}
\]
be the least digit occurring infinitely many times in the standard
decimal expansion of~$\pi$.  Similarly, let
\[
d^\infty_{\min}
\]
be the least digit occurring in an infinite place of the extended
decimal expansion of~$\pi$.  Then the following equality holds:
\[
d_{\min} = d^\infty_{\min} .
\]

This equality indicates that our scant knowledge of the infinite
decimal places of~$\pi$ is not due entirely to the ``non-constructive
nature of the classical constructions using the axiom of choice", as
has sometimes been claimed; but rather to our scant knowledge of the
standard decimal expansion: no ``naturally arising'' irrationals are
known to possess infinitely many occurrences of any specific digit.

\begin{question}%2.13
What does the odd expression ``$H$-infinitely many"
mean exactly?
\end{question}

Answer.  A typical application of an infinite hyperinteger~$H$ is the
proof of the extreme value theorem.%
\footnote{
See Appendix~\ref{ten}, item~\ref{109} for details.}
Here one partitions the interval, say~$[0,1]$, into~$H$-infinitely
many equal subintervals (each subinterval is of course infinitesimally
short).  Then we find the maximum~$x_{i_0}$ among the~$H+1$ partition
points~$x_i$ by the transfer principle,%
\footnote{
See Appendix~\ref{ten}, item~\ref{101}.}
and point out that by continuity, the standard part of the
hyperreal~$x_{i_0}$ gives a maximum of the real function.

\begin{question}%
I am still bothered by changing the meaning of the
notation~$.999\ldots$ as it can be misleading.  I recall I was taught
that it is preferable to use the~$y'$ or~$y_x$ notation until one is
familiar with derivatives, since~$dy/dx$ can be very misleading even
though it can be extremely convenient.  Shouldn't it be avoided?
\end{question}

Answer.  There may be a reason for what you were taught, already noted
by Bishop Berkeley%
\footnote{
See footnote~\ref{rolle} and main text there.}
\cite{Be} nearly 300 years ago!  Namely, standard analysis has no way
of justifying these manipulations rigorously.  The introduction of the
notation~$dy/dx$ is postponed in the standard approach, until the
students are already comfortable with derivatives, as the implied {\em
ratio\/} is thought of as misleading.

Meanwhile, mathematician and leading mathematics educator D.~Tall
writes as follows \cite[chapter~11]{Ta10}:
\begin{quote}
What is far more appropriate for beginning students is an approach
building from experience of dynamic embodiment and the familiar
manipulation of symbols in which the idea of $dy/dx$ as the ratio of
the components of the tangent vector is fully meaningful.
\end{quote}
Tall writes that, following the adoption of the limit concept by
mathematicians as the basic one, 
\begin{quote}
[s]tudents were given emotionally charged instructions to avoid
thinking of $dy/dx$ as a ratio, because it was now seen as a limit,
even though the formulae of the calculus operated as if the expression
were a ratio, and the limit concept was intrinsically problematic.
\end{quote}
With the introduction of infinitesimals such as~$\Delta x$, one
defines the derivative~$f'(x)$ as
\[
f '(x) = st(\Delta y / \Delta x),
\]
where ``st" is the standard part function.%
\footnote{
See Appendix~\ref{ten}, item~\ref{103} and item~\ref{105}.}
Then one sets~$dx=\Delta x$, and defines~$dy=f'(x)dx$.  Then~$f '(x)$
is truly the ratio of two infinitesimals:~$f'(x)=dy/dx$, as envisioned
by the founders of the calculus%
\footnote{
\label{bernoulli}Schubring \cite[p.~170, 173, 187]{Sch}
attributes the first {\em systematic\/} use of infinitesimals as a
foundational concept, to Johann Bernoulli (rather than Newton or
Leibniz).}
and justified by Robinson.

\section{A cardinal issue}%6

\begin{question}%
\label{q61}
How does one relate hyperreal infinities to cardinality?  It still
isn't clear to me what ``$H$-infinitely many 9s" means.  Is
it~$\aleph_0$,~$\aleph_1$, the continuum, or something else?
\end{question}

Answer.  Since there exist {\em countable\/} Skolem non-standard
models of arithmetic \cite{Sk}, the {\em short\/} answer to your
question is ``$\aleph_0$''.  Every non-standard natural number in such
a model will of course have only countably many numbers smaller than
itself, and therefore every extended decimal will have only countably
many digits.

\begin{question}
\label{skolem}
How do you go from Skolem to {\em point, nine recurring\/}?
\end{question}

Answer.  Skolem \cite{Sk} already constructed non-standard models of
arithmetic a quarter century before A. Robinson.  Following the work
of J.~Avigad \cite{Av01}, it is possible to capture a significant
fragment of non-standard calculus, in a very weak logical language;
namely, in the language of primitive recursive arithmetic (PRA), in
the context of the {\em fraction field\/} of Skolem's non-standard
model.  Avigad gives an explicit syntactic translation of the
nonstandard theory to the standard theory.  In the fraction field of
Skolem's non-standard model, equality, thought of as a two-place
relation, is {\em interpreted\/} as the relation $=_*$ of being
infinitely close.  A ``real number'' can be thought of as an
equivalence class relative to such a relation, though the actual
construction of the quotient space (``the continuum'') transcends the
PRA framework.

The integer part (i.e.~floor) function $[x]$ is primitive recursive
due to the existence of the Euclidean algorithm of long division.
Thus, we have~$[m/n] = 0$ if~$m<n$, and similarly~$[m+n/n]= [m/n]+1$.
Furthermore, the digits of a decimal expansion are easily expressed in
terms of the integer part.  Hence the digits are primitive recursive
functions.  Thus the PRA framework is sufficient for dealing with the
issue of {\em point, nine recurring\/}.  Such an approach provides a
common extended decimal ``kernel'' for most theories containing
infinitesimals, not only Robinson's theory.

\begin{question}%3.1b
What's the {\em long\/} answer on cardinalities?
\end{question}

Answer.  On a deeper level, one needs to get away from the naive
cardinals of Cantor's theory,%
\footnote{
which cannot be used in any obvious way as individual numbers in an
extended number system.}
and focus instead on the distinction between a language and a model.
A language (more precisely, a theory in a language, such as first
order logic) is a collection of propositions.  One then interprets
such propositions with respect to a particular model.

A key notion here is that of an internal set.%
\footnote{
See Appendix~\ref{ten}, item~\ref{102}.}
Each set~$S$ of reals has a natural extension~$S^*$ over~$\R^*$, but
also atomic elements of $\R^*$ are considered internal, so the
collection of internal sets is somewhat larger than just the natural
extensions of real sets.

A key observation is that, when the language is being applied to the
non-standard extension, the propositions are being interpreted as
applying only to internal sets, rather than to all sets.

In more detail, there is a certain set-theoretic construction%
\footnote{
See discussion of the ultrapower construction in
Section~\ref{eight}}
of the hyperreal theory~$\R^*$, but the language will be interpreted
as applying only to internal sets and not all set-theoretic subsets
of~$\R^*$.

Such an interpretation is what makes it possible for the transfer
principle to hold, when applied to a theory in first order language.

\begin{question}%
I still have no idea what the extended decimal expansion is.
\end{question}

Answer.  In Robinson's theory, the set of standard natural
numbers~$\N$ is imbedded inside the collection of hyperreal natural
numbers,%
\footnote{
A non-standard model of arithmetic is in fact sufficient for
our purposes; see the answer to Question~\ref{q61}}
denoted~$\N^*$.  The elements of the difference~$\N^* \setminus \N$
are sometimes called (positive) infinite hyperintegers, or
non-standard integers.
%(see photo of Figure~\ref{robinsonfig}).

The standard decimal expansion is thought of as a string of digits
labeled by~$\N$.  Similarly, Lightstone's extended expansion can be
thought of as a string labeled by~$\N^*$.  Thus an extended decimal
expansion for a hyperreal in the unit interval will appear as
\[
a = .a_1a_2a_3\ldots;\ldots a_{H-2} a_{H-1} a_H\ldots
\]
The digits before the semicolon are the ``standard" ones (i.e.~the
digits of $\st(a)$, see Appendix~\ref{ten}).  Given an infinite
hyperinteger~$H$, the string containing~$H$-infinitely many 9s will be
represented by
\[
.999\ldots;\ldots 999
\]
where the last digit 9 appears in position~$H$.  It falls short of 1
by the infinitesimal amount~$1/10^H$.

\begin{question}%
What happens if one decreases $.999\ldots;\ldots 999$ further, by the
same infinitesimal amount~$1/10^H$ ?
\end{question}

Answer.  One obtains the hyperreal number
\[
.999\ldots;\ldots 998,
\]
with digit ``$8$'' appearing at infinite rank~$H$.

\begin{question}%3.4
You mention that the students have not been taught about~$\R$
and~$\lim$ before being introduced to non-terminating decimals.
Perhaps the best solution is to delay the introduction of
non-terminating decimals?  What point is there in seeking the ``right"
approach, if in any case the students will not know what you are
talking about?
\end{question}

Answer.  How would you propose to implement such a scheme?  More
specifically, just how much are we to divulge to the students about
the result of the long division%
\footnote{
This long division is analyzed in the answer to
Question~\ref{q35} below}
of~$1$ by~$3$?

\begin{question}%
Just between the two of us, in the end, there is still no theoretical
explanation for the strict inequality $.999\ldots < 1$, is there?  You
did not disprove the equality~$.999\ldots = 1$.  Are there any
schoolchildren that could understand Lightstone's notation?
\end{question}

Answer. The point is not to teach Lightstone's notation to
schoolchildren, but to broaden their horizons by mentioning the
existence of arithmetic frameworks where their ``hunch"
that~$.999\ldots$ falls short of~$1$, can be justified in a
mathematically sound fashion, consistent with the idea of an
``infinite string of 9s" they are already being told about.  The
underbrace notation
\begin{equation*}
\label{under}
.\underset{H}{\underbrace{999\ldots}}\; = 1- \frac{1}{10^{H}}
\end{equation*}
may be more self-explanatory than Lightstone's semicolon notation; to
emphasize the infinite nature of the non-standard integer~$H$, one
could denote it by the traditional infinity symbol~$\infty$, so as to
obtain a strict inequality%
\footnote{
See answer to Question~\ref{q36} for a more specific choice
of~$H$}
\begin{equation*}
.\underset{\infty}{\underbrace{999\ldots}} < 1,
\end{equation*}
keeping in mind that the left-hand side is an infinite terminating
extended decimal.

\begin{question}%
The multitude of bad teachers will stumble and misrepresent whatever
notation you come up with.  For typesetting purposes, Lighthouse's
notation is more suitable than your underbrace notation.  Isn't an
able mathematician committing a capital sin by promoting a pet
viewpoint, as the cure-all solution to the problems of math education?
\end{question}

Answer.  Your assessment is that the situation is bleak, and the
teachers are weak.  On the other hand, you seem to be making a hidden
assumption that the status-quo cannot be changed in any way.  Without
curing all ills of mathematics education, one can ask what educators
think of a specific proposal addressing a specific minor ill, namely
student frustration with the problem of unital evaluation.

One solution would be to dodge the discussion of it altogether.  In
practice, this is not what is done, but rather the students are indeed
presented with the claim of the evaluation of .999\ldots to 1.  This is
done before they are taught~$\R$ or~$\lim$.  The facts on the ground
are that such teaching is indeed going on, whether in~$12^{th}$ grade
(or even earlier, see~\cite{WAD}) or at the freshman level.

\begin{question}%
Are hyperreals conceptually easier than the common reals?  Will modern
children interpret sensibly ``infinity minus one," say?
\end{question}

Answer.  David Tall, a towering mathematics education figure, has
published the results of an ``interview" with a pre-teen, who quite
naturally developed a number system where 1, 2, 3 can be added to
``infinity" to obtain other, larger, ``infinities".  This indicates
that the idea is not as counterintuitive as it may seem to us, through
the lens of our standard education.

\begin{question}%
\label{q28}
If the great Kronecker could not digest Cantor's infinities, how are
modern children to interpret them?
\end{question}

Answer.  No, schoolchildren should not be taught the arithmetic of the
hyperreals, no more than Cantorian set theory.  On the other hand, the
study by K. Sullivan \cite{Su} in the Chicago area indicates that
students following the non-standard calculus course were better able
to interpret the sense of the mathematical formalism of calculus than
a control group following a standard syllabus.  Sullivan's conclusions
were also noted by Artigue \cite{Ar}, Dauben~\cite{Da96}, and Tall
\cite{Ta80}.  A more recent synthesis of teaching frameworks based on
non-standard calculus was developed by Bernard Hodgson~\cite{Ho} in
'92, and presented at the ICME-7 at Quebec.

Are these students greater than Kronecker?  Certainly not.  On the
other hand, Kronecker's commitment to the ideology of finitism%
\footnote{
\label{semimystical}As a lightning introduction to
Intuitionism, we note that Kronecker rejected actual (completed)
infinity, as did Brouwer, who also rejected the law of excluded middle
(which would probably have been rejected by Kronecker, had it been
crystallized as an explicit concept by logicians in Kronecker's time).
Brouwer developed a theory of the continuum in terms of his ``choice
sequences''.  E.~Bishop's Constructivism rejects both Kronecker's
finitism (Bishop accepts the actual infinity of $\N$) and Brouwer's
theory of the continuum, described as ``semimystical'' by Bishop
\cite[p.~10]{Bi67}.}
was as powerful as most mathematicians' commitment to the standard
reals is, today.

Mathematics education researcher J.~Monaghan, based on field studies,
has reached the following conclusion \cite[p.~248]{Mo}:
\begin{quote}
[...] do infinite numbers of any form exist for young people without
formal mathematical training in the properties of infinite numbers?
The answer is a qualified `yes'.
\end{quote}

\begin{question}%
Isn't the more sophisticated reader going to wonder why Lightstone
stated in \cite{Li} that decimal representation is unique, while you
are making a big fuss over the nonuniqueness of decimal notation and
the strict inequality?
\end{question}

Answer.  Lightstone was referring to the convention of replacing each
terminating decimal, by a tail of 9s.  

Beyond that, it is hard to get into Lightstone's head.  Necessarily
remaining in the domain of speculation, one could mention the
following points.  Mathematicians trained in standard decimal theory
tend to react with bewilderment to any discussion of a strict
inequality~``$.999\ldots$''~$< 1$.  Now Lightstone was interested in
publishing his popular article on infinitesimals, following his
advisor's (Robinson's) approach.  There is more than one person
involved in publishing an article.  Namely, an editor also has a say,
and one of his priorities is defining the level of controversy
acceptable in his periodical.

\begin{question}%
Why didn't Lighstone write down the strict inequality?
\end{question}

Answer.  Lightstone could have made the point that all but one
extended expansions starting with~$999\ldots$ give a hyperreal value
strictly less than~$1$.  Instead, he explicitly reproduces only the
expansion equal to~$1$.  In addition, he explicitly mentions an
additional expansion--and explains why it does not exist!  Perhaps he
wanted to stay away from the strict inequality, and concentrate
instead on getting a minimal amount of material on non-standard
analysis published in a mainstream popular periodical.  All this is in
the domain of speculation.

As far as the reasons for elaborating a strict non-standard
inequality, they are more specific.  First, the manner in which the
issue is currently handled by education professionals, tends to
engender student frustration.  Furthermore, the standard treatment
conceals the power of non-standard analysis in this particular issue.

\section{Circular reasoning, ultrafilters, and Platonism}%7
\label{eight}

\begin{question}%
\label{q51}
Since the construction of the hyperreal numbers depends on that of the
real numbers, wouldn't it be extremely easy for people to attack this
idea as being circular reasoning?
\end{question} 

Answer.  Actually, your assumption is incorrect.  Just as the reals
can be obtained from the rationals as the set of equivalence classes
of suitable sequences of rational numbers (namely, the Cauchy ones),
so also a version of the hyperreals%
\footnote{
One does not obtain all elements of $\R^*$ by starting with
sequences of rational numbers, but the resulting non-Archimedean
extension of $\R$ is sufficient for most purposes of the calculus,
cf.~Avigad~\cite{Av01}.}
can be obtained from the rationals as the set of equivalence classes
of sequences of rational numbers, modulo a suitable equivalence
relation.  Such a construction is due to Luxemburg~\cite{Lu62}.  The
construction is referred to as the ultrapower construction, see
Goldblatt~\cite{Go}.

\begin{question}%
Non-standard analysis?  You mean ultrafilters and all that?
\end{question}

Answer.  The good news is that you don't need ultrafilters to do
non-standard analysis: the axiom of choice is enough.%
\footnote{
The comment is, of course, tongue in cheek, but many people seem not
to have realized yet that the existence of a free (non-principal)
ultrafilter is as much of a consequence of the axiom of choice, as the
existence of a maximal ideal (a standard tool in algebra), or the
Hahn-Banach theorem (a standard tool in functional analysis).  This is
as good a place as any to provide a brief unequal time to an opposing
view \cite{Ri94}: {\em We are all Platonists, aren't we? In the
trenches, I mean---when the chips are down. Yes, Virginia, there
really are circles, triangles, numbers, continuous functions, and all
the rest.  Well, maybe not free ultrafilters. Is it important to
believe in the existence of free ultrafilters? Surely that's not
required of a Platonist. I can more easily imagine it as a test of
sanity: `He believes in free ultrafilters, but he seems harmless'.}
Needless to say, the author of \cite{Ri94} is in favor of eliminating
the axiom of choice--including the countable one.}

\begin{question}%
\label{q33}
Good, because, otherwise, aren't you sweeping a lot under the rug when
you teach non-standard analysis to first year students?
\end{question}

Answer.  One sweeps no more under the rug than the equivalence classes
of Cauchy sequences, which are similarly {\em not\/} taught in first
year calculus.  After all, the hyperreals are just equivalence classes
of more general sequences
%
%of rational numbers: this was an error
%
(this is known as the ultrapower construction).  What one does {\em
not\/} sweep under the rug in the hyperreal approach is the notion of
infinitesimal which historically was present at the inception of the
theory, whether by Archimedes or Leibniz-Newton.  Infinitesimals were
routinely used in teaching until as late as 1912, the year of the last
edition of the textbook by L. Kiepert~\cite{Ki}.  This issue was
discussed in more detail by P.~Roquette~\cite{Roq}.

\begin{question}%
I have a serious problem with Lightstone's notation.  I can see it
working for a specific infinite integer~$H$, and even for nearby
infinite integers of the form~$H+n$, where~$n$ is a finite integer,
positive or negative.  However, I do not see how it represents two
different integers, for instance~$H$ and~$H^2$ on the same picture.
For in this case,~$H^2$ is greater than~$H+n$ for any finite~$n$.
Thus it does not lie in the same infinite collection of decimal places
\[
;\ldots 1 \ldots
\]
so that one needs even more than a potentially infinite collection of
sequences of digits
\[
;\ldots ;
\]
to cope with all hyperintegers.
\end{question}

Answer.  Skolem \cite{Sk} already constructed non-standard models for
arithmetic, many years before Robinson.%
\footnote{
See the answer to Question~\ref{skolem} for more details.}
Here you have a copy of the standard integers, and also many
``galaxies".  A galaxy, in the context of Robinson's hyperintegers,%
\footnote{
See Goldblatt~\cite{Go} for more details}
is a collection of hyperintegers differing by a finite integer.  At
any rate, one does need infinitely many semicolons if one were to dot
all the i's.

Lightstone is careful in his article to discuss this issue.  Namely,
what is going on to the right of his semicolon is not similar to the
simple picture to the left.  At any rate, the importance of his
article is that he points out that there does exist the notion of an
extended decimal representation, where the leftmost galaxy of digits
of $x$ are the usual finite digits of $\st(x)$.

\section{How long a division?}%8

\begin{question}%
\label{q35}
Long division of $1$ by $3$ gives $.333\ldots$ which is a very {\bf
obvious pattern}.  Therefore multiplying back by $3$ we get
$.999\ldots=1$.  There is nothing else to discuss!
\end{question}

Answer.  Let us be clear about one thing: long division of~$1$ by~$3$
does {\em not\/} produce the infinite decimal~$.333\ldots$ contrary to
popular belief.  What it {\em does\/} produce is the sequence~$\langle
.3, .33, .333, \ldots \rangle$, where the dots indicate the {\em
obvious pattern\/}.

Passing from a sequence to an infinite decimal is a major additional
step.  The existence of an infinite decimal expansion is a non-trivial
matter that involves the construction of the real number system, and
the notion of the limit.

\begin{question}%
\label{q45}
Doesn't the standard formula for converting every repeated decimal to
a fraction show that $.333\ldots$ equals $\frac{1}{3}$ on the nose?
\end{question}

Answer.  Converting decimals to fractions was indeed the approach
of~\cite{WAD}.  However, in a pre-$\R$ environment, one can argue that
the formula only holds up to an infinitesimal error, and attempts to
``prove'' unital evaluation by an appeal to such a formula amount to
replacing one article of faith, by another.

To elaborate, note that applying the iterative procedure of long
division in the case of~$\frac{1}{3}$, does not by itself produce any
infinite decimal, no more than the iterative procedure of adding~$1$
to the outcome of the previous step, produces any infinite integer.
Rather, the long division produces the sequence~$\langle .3, .33,
.333, \ldots \rangle$.  Transforming the sequence into an infinite
decimal has nothing to do with long division, and requires, rather, an
application of the limit concept, in the context of a complete number
system.

Note that, if we consider the sequence
\[
\langle .9, .99, .999, \ldots \rangle ,
\]
but instead of taking the limit, take its equivalence class
\[
[.9, .99, .999, \ldots]
\]
in the ultrapower construction of the hyperreals (see \cite{Lu62,
Go}), then we obtain a value equal to~$1 - 1/10^{[\N]}$ where~$\langle
\N \rangle$ is the ``natural string" sequence enumerating all the {\em
natural\/} numbers, whereas~$[\N]$ is its equivalence class in the
hyperreals.  Thus the unital evaluation has a viable competitor,
namely, the ``natural string'' evaluation.

\begin{question}%
\label{q36}
Isn't it odd that you seem to get a {\em canonical\/} representative
for ``$.999\ldots$'' which falls short of $1$?
\end{question}

Answer.  The hyperinteger defined by the equivalence class of the
sequence $\langle\N\rangle=\langle 1,2,3,\ldots\rangle$ only makes
sense in the context of the ultrapower construction, and depends on
the choices made in the construction.  The standard real
decimal~$(.999\ldots)_{\rm Lim}$ is defined as the {\em limit\/} of
the sequence~$(.9, .99, .999 \ldots)$, and the
hyperreal~$(.999\ldots)_{\rm Lux}$ is defined as the class of the same
sequence in the ultrapower construction.  Then
\begin{equation}
(.999\ldots)_{\rm Lim}=1
\end{equation}
is the unital evaluation, interpreting the symbol~$.999\ldots$ as a
real number, while
\begin{equation}
.\underset{[\N]}{\underbrace{999\ldots}} = 1 - \frac{1}{10^{[\N]}}.
\end{equation}
is the natural string evaluation.

In more detail, in the ultrapower construction of~$\R^*$, the
hyperreal~$[.9, .99, .999, \ldots]$, represented by the
sequence~$\langle .9, .99, .999, \ldots \rangle$, is an infinite
terminating string of~$9$s, with the last nonzero digit occurring at a
suitable infinite hyperinteger rank.  The latter is represented by the
string listing all the natural numbers~$\langle 1,2,3,\ldots \rangle$,
which we abbreviate by the symbol~$\langle \N \rangle$.  Then the
equivalence class~$[\N]$ is the corresponding ``natural string''
hyperinteger.  We therefore obtain a hyperreal equal
to~$1-\frac{1}{10^{[\N]}}$.

The {\em unital evaluation\/} of the symbol~$.999\ldots$ has a viable
competitor, namely the {\em natural string evaluation\/} of this
symbol.

\begin{question}%
\label{q46}
What do you mean by~$10^{\N}$?  It looks to me like a typical
sophomoric error.
\end{question}

Answer.  The natural string evaluation yields a hyperreal with
Lightstone~\cite{Li} representation given by~$.999\ldots;\ldots 9$,
with the last digit occurring at non-standard rank~$[\N]$.  Note that
it would be {\em incorrect\/} to write
\begin{equation}
(.999\ldots)_{\rm Lux}= 1 - \frac{1}{10^{\N}}
\end{equation}
since the expression~$10^{\N}$ is meaningless,~$\N$ not being a number
in {\em any\/} number system.  Meanwhile, the sequence~$ \langle \N
\rangle$ listing all the natural numbers in increasing order,
represents an equivalence class~$[\N]$ in the ultrapower construction
of the hyperreals, so that $[\N]$ is indeed a quantity, more precisely
a non-standard integer, or a hypernatural number \cite{Go}.

\begin{question}%
Do the subscripts in~$(.999\ldots)_{\rm Lim}$ and~$(.999\ldots)_{\rm
Lux}$ stand for ``limited'' and ``deluxe''?
\end{question}

Answer.  No, the subscript ``Lim'' refers to the unital evaluation
obtained by applying the {\em limit\/} to the sequence, whereas the
subscript ``Lux'' refers to the natural string evaluation, in the
context of Luxemburg's sequential construction of the hyperreals (the
ultrapower construction).

\begin{question}%
The absence of infinitesimals is certainly not some kind of a
shortcoming of the real number system that one would need to apologize
for.  How can you imply otherwise?
\end{question}

Answer. The standard reals are at the foundation of the magnificent
edifice of classical and modern analysis.  Ever since their rigorous
conception by Weierstrass, Dedekind, and Cantor, the standard reals
have faithfully served the needs of generations of mathematicians of
many different specialties.  Yet the non-availability of
infinitesimals has the following consequences:
\begin{enumerate}
\item
it distances mathematics from its applications in physics, enginering,
and other fields (where nonrigorous infinitesimals are in routine
use);
\item
it complicates the logical structure of calculus concepts (such as the
limit) beyond the comprehension of a significant minority (if not a
majority) of undergraduate students;
\item
it deprives us of a key tool in interpreting the work of such greats
as Archimedes, Euler, and Cauchy. 
\end{enumerate}
In this sense, the absence of infinitesimals is a shortcoming of the
standard number system.

\appendix

\section{A non-standard glossary}
\label{ten}

The present section can be retained or deleted at the discretion of
the referee.  In this section we present some illustrative terms and
facts from non-standard calculus \cite{Ke}.  The relation of being
infinitely close is denoted by the symbol~$\approx$.  Thus,~$x\approx
y$ if and only if~$x-y$ is infinitesimal.

\subsection{Natural hyperreal extension~$f^*$}
\label{101}
The {\em extension principle\/} of non-standard calculus states that
every real function~$f$ has a hyperreal extension, denoted~$f^*$ and
called the natural extension of~$f$.  The {\em transfer principle\/}
of non-standard calculus asserts that every real statement true
for~$f$, is true also for~$f^*$ (for statements involving any
relations).  For example, if~$f(x)>0$ for every real~$x$ in its
domain~$I$, then~$f^*(x)>0$ for every hyperreal~$x$ in its
domain~$I^*$.  Note that if the interval~$I$ is unbounded, then~$I^*$
necessarily contains infinite hyperreals.  We will sometimes drop the
star~$^*$ so as not to overburden the notation.

\subsection{Internal set}
\label{102}
Internal set is the key tool in formulating the transfer principle,
which concerns the logical relation between the properties of the real
numbers~$\R$, and the properties of a larger field denoted
\[
\R^*
\]
called the {\em hyperreal line}.  The field~$\R^*$ includes, in
particular, infinitesimal (``infinitely small") numbers, providing a
rigorous mathematical realisation of a project initiated by Leibniz.
Roughly speaking, the idea is to express analysis over~$\R$ in a
suitable language of mathematical logic, and then point out that this
language applies equally well to~$\R^*$.  This turns out to be
possible because at the set-theoretic level, the propositions in such
a language are interpreted to apply only to internal sets rather than
to all sets.  Note that the term ``language" is used in a loose sense
in the above.  A more precise term is {\em theory in first-order
logic}.  Internal sets include natural extension of standard sets.

\subsection{Standard part function}
\label{103}
The standard part function ``st" is the key ingredient in
A.~Robinson's resolution of the paradox of Leibniz's definition of the
derivative as the ratio of two infinitesimals
\[
\frac{dy}{dx}.
\]
The standard part function associates to a finite hyperreal
number~$x$, the standard real~$x_0$ infinitely close to it, so that we
can write
\begin{equation*}
\mathrm{st}(x)=x_0.
\end{equation*}
In other words, ``st'' strips away the infinitesimal part to produce
the standard real in the cluster (see item below).  The standard part
function ``st" is not defined by an internal set (see item~\ref{102}
above) in Robinson's theory.

\subsection{Cluster}
\label{104}
Each standard real is accompanied by a cluster%
\footnote{Alternative terms are prevalent in the literature, such as
{\em halo\/} \cite{Go}, but the term {\em cluster\/} has the advantage
of being self-explanatory}
of hyperreals infinitely close to it.  The standard part function
collapses the entire cluster back to the standard real contained in
it.  The cluster of the real number~$0$ consists precisely of all the
infinitesimals.  Every infinite hyperreal decomposes as a triple sum
\[
H+r+\epsilon,
\]
where~$H$ is a hyperinteger,~$r$ is a real number in~$[0,1)$, and
$\epsilon$ is infinitesimal.  Varying~$\epsilon$ over all
infinitesimals, one obtains the cluster of~$H+r$.

\subsection{Derivative}
\label{105}
To define the real derivative of a real function~$f$ in this approach,
one no longer needs an infinite limiting process as in standard
calculus.  Instead, one sets
\begin{equation}
\label{deri}
f'(x) = \mathrm{st} \left( \frac{f(x+\epsilon)-f(x)}{\epsilon}
\right),
\end{equation}
where~$\epsilon$ is infinitesimal, yielding the standard real number
in the cluster of the hyperreal argument of ``st'' (the derivative
exists if and only if the value~\eqref{deri} is independent of the
choice of the infinitesimal).  The addition of ``st'' to the formula
resolves the centuries-old paradox famously criticized by George
Berkeley%
\footnote{
See footnote~\ref{rolle} for a historical clarification}
\cite{Be} (in terms of the {\em Ghosts of departed
quantities}, cf.~\cite[Chapter~6]{St}), and provides a rigorous basis
for infinitesimal calculus as envisioned by Leibniz.

\subsection{Continuity}
\label{106}
A function~$f$ is continuous at~$x$ if the following condition is
satisfied:~$y\approx x$ implies~$f(y)\approx f(x)$.

\subsection{Uniform continuity} 
\label{107}
A function~$f$ is uniformly continuous on~$I$ if the following
condition is satisfied:

\begin{itemize}
\item
standard: for every~$\epsilon>0$ there exists a~$\delta>0$ such that
for all~$x\in I$ and for all~$y\in I$, if~$|x-y|<\delta$ then
$|f(x)-f(y)| < \epsilon$.
\item
non-standard: for all~$x\in I^*$, if~$x\approx y$ then~$f(x) \approx
f(y)$.
\end{itemize}
%See Remark~\ref{51b} for a more detailed discussion.

\subsection{Hyperinteger}
\label{108}
A hyperreal number~$H$ equal to its own integer part 
\[
H = [H]
\]
is called a hyperinteger (here the integer part function is the
natural extension of the real one).  The elements of the complement
$\Z^* \setminus \Z$ are called infinite hyperintegers, or non-standard
integers.
%, see Figure~\ref{robinsonfig}.

\subsection{Proof of extreme value theorem}
\label{109}
Let~$H$ be an infinite hyperinteger.  The interval~$[0,1]$ has a
natural hyperreal extension.  Consider its partition into~$H$
subintervals of equal length~$\frac{1}{H}$, with partition points~$x_i
= i/H$ as~$i$ runs from~$0$ to~$H$.  Note that in the standard
setting, with~$n$ in place of~$H$, a point with the maximal value
of~$f$ can always be chosen among the~$n+1$ partition points~$x_i$, by
induction.  Hence, by the transfer principle, there is a
hyperinteger~$i_0$ such that~$0\leq i_0 \leq H$ and
\begin{equation}
\label{101b}
f(x_{i_0})\geq f(x_i) \quad \forall i= 0,\ldots,H.
\end{equation}
Consider the real point
\begin{equation*}
c= {\rm st}(x_{i_0}).
\end{equation*}
An arbitrary real point~$x$ lies in a suitable sub-interval of the
partition, namely~$x\in [x_{i-1},x_i]$, so that~${\rm st}(x_i) = x$.
Applying ``st'' to the inequality \eqref{101b}, we obtain by continuity
of~$f$ that~$f(c)\geq f(x)$, for all real~$x$, proving~$c$ to be a
maximum of~$f$ (see \cite[p.~164]{Ke}).

\subsection{Limit}
\label{1010}
We have~$\lim_{x\to a} f(x) = L$ if and only if whenever the
difference~$x-a\not=0$ is infinitesimal, the difference~$f(x)-L$ is
infinitesimal, as well, or in formulas: if~${\rm st}(x)=a$ then~${\rm
st}(f(x)) = L$.

Given a sequence of real numbers~$\{x_n|n\in \mathbb{N}\}$, if~$L\in
\mathbb{R}\;$ we say~$L$ is the limit of the sequence and write~$L =
\lim_{n \to \infty} x_n$ if the following condition is satisfied:
\begin{equation}
\label{102b}
{\rm st} (x_H)=L \quad \mbox{\rm for all infinite } H
\end{equation}
(here the extension principle is used to define~$x_n$ for every
infinite value of the index).  This definition has no quantifier
alternations.  The standard~$(\epsilon, \delta)$-definition of limit,
on the other hand, does have quantifier alternations:
\begin{equation}
\label{disaster}
L = \lim_{n \to \infty} x_n\Longleftrightarrow \forall \epsilon>0\;,
\exists N \in \mathbb{N}\;, \forall n \in \mathbb{N} : n >N \implies
d(x_n,L)<\epsilon.
\end{equation}          

\subsection{Non-terminating decimals}
\label{1011}
Given a real decimal $u=.d_1d_2d_3\ldots$, consider the
sequence~$u_1=.d_1$,~$\;u_2=.d_1 d_2$,~$\;u_3=.d_1 d_2 d_3$, etc.
Then by definition,
\[
u=\lim_{n\to \infty} u_n.
\]
Meanwhile,~$\lim_{n\to \infty} u_n= \st(u_H^{\phantom{I}})$ for every
infinite~$H$.  Now if~$u$ is a non-terminating decimal, then one
obtains a strict inequality~$u_H^{\phantom{I}}<u$ by transfer from
$u_n<u$.  In particular,
\begin{equation}
\label{105b}
.999\ldots;\ldots \hat 9 = .\underset{H}{\underbrace{999\ldots}} = 1 -
\tfrac{1}{10^H} < 1,
\end{equation}
where the hat~$\hat {\phantom{9}}$ indicates the~$H$-th Lightstone
decimal place.  The standard interpretation of the symbol~$.999\ldots$
as~$1$ is necessitated by notational uniformity: the symbol~$.a_1 a_2
a_3 \ldots$ in every case corresponds to the limit of the sequence of
terminating decimals~$.a_1\ldots a_n$.  Alternatively, the ellipsis
in~$.999\ldots$ could be interpreted as alluding to an infinity of
nonzero digits specified by a choice of an infinite hyperinteger~$H\in
\N^*\setminus \N$.  The resulting~$H$-infinite extended decimal string
of~$9$s corresponds to an infinitesimally diminished hyperreal
value~\eqref{105b}.  Such an interpretation is perhaps more in line
with the naive initial intuition persistently reported by teachers.

\subsection{Integral} 
\label{1012}
The definite integral of~$f$ is the standard part of an infinite
Riemann sum~$\sum_{i=0}^H f(x) \Delta x$, the latter being defined by
means of the transfer principle, once finite Riemann sums are in
place, see \cite{Ke} for details.

\section*{Acknowledgments}

We are grateful to Alexander Bogomolny, Norman Megill, and Fred
Richman for many probing questions that have made this text possible.

\end{document}